\input amstex
\documentstyle{amsppt}
%
\catcode`@=11
\redefine\output@{%
  \def\break{\penalty-\@M}\let\par\endgraf
  \ifodd\pageno\global\hoffset=105pt\else\global\hoffset=8pt\fi  
  \shipout\vbox{%
    \ifplain@
      \let\makeheadline\relax \let\makefootline\relax
    \else
      \iffirstpage@ \global\firstpage@false
        \let\rightheadline\frheadline
        \let\leftheadline\flheadline
      \else
        \ifrunheads@ 
        \else \let\makeheadline\relax
        \fi
      \fi
    \fi
    \makeheadline \pagebody \makefootline}%
  \advancepageno \ifnum\outputpenalty>-\@MM\else\dosupereject\fi
}
\def\Beta{\mathchar"0\hexnumber@\rmfam 42}
\catcode`\@=\active
\nopagenumbers
\chardef\textvolna='176

\chardef\bigalpha='013
\def\negskp{\hskip -2pt}
\def\Rea{\operatorname{Re}}
\def\Img{\operatorname{Im}}

\chardef\degree="5E
\def\blue#1{#1}

\gdef\darkred#1{#1}
\catcode`#=11\def\diez{#}\catcode`#=6
\catcode`&=11\catcode`&=4
\catcode`_=11\def\podcherkivanie{_}\catcode`_=8
\catcode`\^=11\catcode`\^=7
\catcode`~=11\def\volna{~}\catcode`~=\active
\def\mycite#1{\cite{\blue{#1}}\immediate\special{ps:
     ShrHPSdict begin /ShrBORDERthickness 0 def}}
\def\myciterange#1#2#3#4{\cite{\blue{#2#3#4}}\immediate\special{ps:
     ShrHPSdict begin /ShrBORDERthickness 0 def}}
\def\mytag#1{%
    \tag#1}
\def\mythetag#1{\thetag{\blue{#1}}\immediate\special{ps:
     ShrHPSdict begin /ShrBORDERthickness 0 def}}
\def\myrefno#1{\no#1}
\def\myhref#1#2{\blue{#2}\immediate\special{ps:
     ShrHPSdict begin /ShrBORDERthickness 0 def}}
\def\myEarXivlink{\myhref{http://arXiv.org}{http:/\negskp/arXiv.org}}

\def\mytheorem#1{\csname proclaim\endcsname{Theorem #1}}
\def\mytheoremwithtitle#1#2{\csname proclaim\endcsname{Theorem #1#2}}
\def\mythetheorem#1{\blue{#1}\immediate\special{ps:
     ShrHPSdict begin /ShrBORDERthickness 0 def}}
\def\mylemma#1{\csname proclaim\endcsname{Lemma #1}}
\def\mylemmawithtitle#1#2{\csname proclaim\endcsname{Lemma #1#2}}
\def\mythelemma#1{\blue{#1}\immediate\special{ps:
     ShrHPSdict begin /ShrBORDERthickness 0 def}}
\def\mycorollary#1{\csname proclaim\endcsname{Corollary #1}}

\def\myconjecture#1{\csname proclaim\endcsname{Conjecture #1}}
\def\myconjecturewithtitle#1#2{\csname proclaim\endcsname{Conjecture #1#2}}
\def\mytheconjecture#1{\blue{#1}\immediate\special{ps:
     ShrHPSdict begin /ShrBORDERthickness 0 def}}
\def\myproblem#1{\csname proclaim\endcsname{Problem #1}}
\def\myproblemwithtitle#1#2{\csname proclaim\endcsname{Problem #1#2}}


\pagewidth{360pt}
\pageheight{606pt}
\topmatter
\title
Reverse asymptotic estimates for roots of the cuboid characteristic equation
in the case of the second cuboid conjecture.
\endtitle
\rightheadtext{Reverse asymptotic estimates \dots}
\author
Ruslan Sharipov
\endauthor
\address Bashkir State University, 32 Zaki Validi street, 450074 Ufa, Russia
\endaddress
\email
\myhref{mailto:r-sharipov\@mail.ru}{r-sharipov\@mail.ru}
\endemail
\abstract
     A perfect cuboid is a rectangular parallelepiped whose edges, whose face 
diagonals, and whose space diagonal are of integer lengths. The second cuboid 
conjecture specifies a subclass of perfect cuboids described by one Diophantine 
equation of tenth degree and claims their non-existence within this subclass. 
This Diophantine equation has two parameters. Previously asymptotic expansions
and estimates for roots of this equation were obtained in the case where the 
first parameter is fixed and the other tends to infinity. In the present paper 
reverse asymptotic expansions and estimates are derived in the case where the 
second parameter is fixed and the first one tends to infinity. Their application
to the perfect cuboid problem is discussed.
\endabstract
\endtopmatter
\TagsOnRight
\document

\head
1. Introduction.
\endhead
     Perfect cuboids are described by six Diophantine equations. These
equations are elementary. They are derived with the use of the Pythagorean 
theorem:  
$$
\xalignat 2
&\hskip -2em
x_1^2+x_2^2+x_3^2=L^2,
&&x_2^2+x_3^2=d_1^{\kern 1pt 2},\\
\vspace{-1.7ex}
\mytag{1.1}\\
\vspace{-1.7ex}
&\hskip -2em
x_3^2+x_1^2=d_2^{\kern 1pt 2},
&&x_1^2+x_2^2=d_3^{\kern 1pt 2}.
\endxalignat
$$
Here $x_1$, $x_2$, $x_3$ are the lengths of the edges of a cuboid, 
$d_1$, $d_2$, $d_3$ are the lengths of its face diagonals and $L$ 
is the length of its space diagonal.\par
     None of the solutions for the equations \mythetag{1.1} are
known. Their non-existence is also not proved. This is an open
mathematical problem. For the history and various approaches to 
the problem of perfect cuboids the reader is referred to 
\myciterange{1}{1}{--}{43}. In the present paper we continue the 
research from \myciterange{44}{44}{--}{49}. The series of papers 
\myciterange{50}{50}{--}{62} is devoted to another approach. This
approach is not considered here.\par
      In \mycite{44} an algebraic parametrization for the Diophantine equations 
\mythetag{1.1} was suggested. It uses four rational variables 
$\alpha$, $\beta$, $\upsilon$, and $z$:
$$
\xalignat 2
&\hskip -2em 
\frac{x_1}{L}=\frac{2\,\upsilon}{1+\upsilon^2},
&&\frac{d_1}{L}=\frac{1-\upsilon^2}{1+\upsilon^2},\\
\vspace{1ex}
&\hskip -2em
\frac{x_2}{L}=\frac{2\,z\,(1-\upsilon^2)}{(1+\upsilon^2)\,(1+z^2)},
&&\frac{x_3}{L}=\frac{(1-\upsilon^2)\,(1-z^2)}{(1+\upsilon^2)\,(1+z^2)},
\mytag{1.2}\\
\vspace{1ex}
&\hskip -2em
\frac{d_2}{L}=\frac{(1+\upsilon^2)\,(1+z^2)+2\,z(1-\upsilon^2)}
{(1+\upsilon^2)\,(1+z^2)}\,\beta,
&&\frac{d_3}{L}=\frac{2\,(\upsilon^2\,z^2+1)}{(1+\upsilon^2)\,(1+z^2)}\,\alpha.
\quad
\endxalignat
$$
The variables $\alpha$, $\beta$, $\upsilon$, and $z$ in \mythetag{1.2} are not
independent. The variable $\upsilon$ is expressed through $\alpha$ and $\beta$ 
as a solution of the following algebraic equation:
$$
\gathered
\upsilon^4\,\alpha^4\,\beta^4+(6\,\alpha^4\,\upsilon^2\,\beta^4
-2\,\upsilon^4\,\alpha^4\,\beta^2-2\,\upsilon^4\,\alpha^2
\,\beta^4)+(4\,\upsilon^2\,\beta^4\,\alpha^2+\\
+\,4\,\alpha^4\,\upsilon^2\,\beta^2-12\,\upsilon^4\,\alpha^2\,\beta^2
+\upsilon^4\,\alpha^4+\upsilon^4\,\beta^4
+\alpha^4\,\beta^4)+(6\,\alpha^4\,\upsilon^2+6\,\upsilon^2\,\beta^4-\\
-\,8\,\alpha^2\,\beta^2\,\upsilon^2-2\,\upsilon^4\,\alpha^2
-2\,\upsilon^4\,\beta^2-2\,\alpha^4\,\beta^2-2\,\beta^4\,\alpha^2)
+(\upsilon^4+\beta^4+\\
+\,\alpha^4+4\,\alpha^2\,\upsilon^2+4\,\beta^2\,\upsilon^2
-12\,\beta^2\,\alpha^2)+(6\,\upsilon^2-2\,\alpha^2-2\,\beta^2)+1=0.
\endgathered
\quad
\mytag{1.3}
$$
Then the variable $z$ is expressed through $\alpha$, $\beta$, and $\upsilon$
by the formula 
$$
\hskip -2em
z=\frac{(1+\upsilon^2)\,(1-\beta^2)\,(1+\alpha^2)}{2\,(1+\beta^2)\,(1
-\alpha^2\,\upsilon^2)}.
\mytag{1.4}
$$
The equation \mythetag{1.3} and the formula \mythetag{1.4} mean that we 
have two algebraic functions
$$
\xalignat 2
&\hskip -2em
\upsilon=\upsilon(\alpha,\beta),
&&z=z(\alpha,\beta).
\mytag{1.5}
\endxalignat
$$
Substituting \mythetag{1.5} into \mythetag{1.2}, we get six algebraic
functions
$$
\xalignat 3
&\hskip -2em
x_1=x_1(\alpha,\beta,L),
&&x_2=x_2(\alpha,\beta,L),
&&x_3=x_3(\alpha,\beta,L),
\qquad\\
\vspace{-1.5ex}
\mytag{1.6}\\
\vspace{-1.5ex}
&\hskip -2em
d_1=d_1(\alpha,\beta,L),
&&d_2=d_2(\alpha,\beta,L),
&&d_3=d_3(\alpha,\beta,L),
\qquad
\endxalignat
$$
They are linear with respect to $L$. The functions \mythetag{1.6} satisfy 
the cuboid equations \mythetag{1.1} identically with respect to $\alpha$, 
$\beta$, and $L$ (see Theorem~5.2 in \mycite{44}).\par
     The rational numbers $\alpha$, $\beta$, and $\upsilon$ cam be brought 
to a common denominator:
$$
\xalignat 3
&\hskip -2em
\alpha=\frac{a}{t},
&&\beta=\frac{b}{t},
&&\upsilon=\frac{u}{t}.
\mytag{1.7}
\endxalignat
$$
Substituting \mythetag{1.7} into \mythetag{1.3}, one easily derives the
Diophantine equation
$$
\gathered
t^{12}+(6\,u^2\,-2\,a^2\,-2\,b^2)\,t^{10}
+(u^4\,+b^4+a^4+4\,a^2\,u^2+\\
+\,4\,b^2\,u^2-12\,b^2\,a^2)\,t^8+(6\,a^4\,u^2+6\,u^2\,b^4-8\,a^2\,b^2\,u^2-\\
-\,2\,u^4\,a^2-2\,u^4\,b^2-2\,a^4\,b^2-2\,b^4\,a^2)\,t^6+(4\,u^2\,b^4\,a^2+\\
+\,4\,a^4\,u^2\,b^2-12\,u^4\,a^2\,b^2+u^4\,a^4+u^4\,b^4+a^4\,b^4)\,t^4+\\
+\,(6\,a^4\,u^2\,b^4-2\,u^4\,a^4\,b^2-2\,u^4\,a^2
\,b^4)\,t^2+u^4\,a^4\,b^4=0.
\endgathered
\quad
\mytag{1.8}
$$
The following theorem claims exact relation of the equation \mythetag{1.8} to
the cuboid equations \mythetag{1.1} (see Theorem~5.2 in \mycite{44} and
Theorem 4.1 in \mycite{45}). 
\mytheorem{1.1} A perfect cuboid does exist if and only if for some
positive coprime integer numbers $a$, $b$, and $u$ the Diophantine
equation \mythetag{1.8} has a positive solution $t$ obeying the inequalities
$t>a$, $t>b$, $t>u$, and $(a+t)\,(b+t)>2\,t^2$. 
\endproclaim
     Due to Theorem~\mythetheorem{1.1} it is natural to call \mythetag{1.8}
the cuboid characteristic equation, though this term was not previously 
used for this equation.\par
     In the case of the second cuboid conjecture the parameters $a$, $b$,
and $u$ are related to each other according to one of the two formulas:
$$
\xalignat 2
&\hskip -2em
b\,u=a^2,
&&a\,u=b^2.
\mytag{1.9}
\endxalignat
$$
The first equality \mythetag{1.9} is resolved by means of the formulas
$$
\xalignat 3
&\hskip -2em
a=p\,q,
&&b=p^{\kern 1pt 2}, &&u=q^{\kern 0.7pt 2}.
\mytag{1.10}
\endxalignat
$$
Here $p\neq q$ are two positive coprime integers. Upon substituting 
\mythetag{1.10} into the equation \mythetag{1.8} it reduces to the 
equation 
$$
(t-a)\,(t+a)\,Q_{p\kern 0.6pt q}(t)=0
\mytag{1.11}
$$
(see \mycite{45} or \mycite{47}), where $Q_{p\kern 0.6pt q}(t)$ is the following 
polynomial of tenth degree:
$$
\gathered
Q_{p\kern 0.6pt q}(t)=t^{10}+(2\,q^{\kern 0.7pt 2}+p^{\kern 1pt 2})\,(3
\,q^{\kern 0.7pt 2}-2\,p^{\kern 1pt 2})\,t^8+(q^{\kern 0.5pt 8}+10
\,p^{\kern 1pt 2}\,q^{\kern 0.5pt 6}+\\
+\,4\,p^{\kern 1pt 4}\,q^{\kern 0.5pt 4}-14\,p^{\kern 1pt 6}\,q^{\kern 0.7pt 2}
+p^{\kern 1pt 8})\,t^6-p^{\kern 1pt 2}\,q^{\kern 0.7pt 2}
\,(q^{\kern 0.5pt 8}-14\,p^{\kern 1pt 2}\,q^{\kern 0.5pt 6}
+4\,p^{\kern 1pt 4}\,q^{\kern 0.5pt 4}+\\
+\,10\,p^{\kern 1pt 6}\,q^{\kern 0.7pt 2}+p^{\kern 1pt 8})\,t^4
-p^{\kern 1pt 6}\,q^{\kern 0.5pt 6}\,(q^{\kern 0.7pt 2}
+2\,p^{\kern 1pt 2})\,(3\,p^{\kern 1pt 2}-2\,q^{\kern 0.7pt 2})\,t^2
-q^{\kern 0.7pt 10}\,p^{\kern 1pt 10}.
\endgathered\quad
\mytag{1.12}
$$\par
The second equality \mythetag{1.9} is similar yo the first one. It is resolved 
by means of the following formulas similar to \mythetag{1.10}: 
$$
\xalignat 3
&\hskip -2em
a=p^{\kern 1pt 2},
&&b=p\,q, &&u=q^{\kern 0.7pt 2}.
\mytag{1.13}
\endxalignat
$$
Upon substituting \mythetag{1.13} into the equation \mythetag{1.8} it reduces 
to the equation 
$$
\hskip -2em
(t-b)\,(t+b)\,Q_{p\kern 0.6pt q}(t)=0.
\mytag{1.14}
$$
The roots $t=a$, $t=-a$, $t=b$, and $t=-b$ of the equations \mythetag{1.11}
and \mythetag{1.14} do not produce perfect cuboids 
(see Theorem~\mythetheorem{1.1}). Upon splitting off the inessential linear 
factors from \mythetag{1.11} and \mythetag{1.14} we get the equation 
$$
\gathered
t^{10}+(2\,q^{\kern 0.7pt 2}+p^{\kern 1pt 2})\,(3\,q^{\kern 0.7pt 2}
-2\,p^{\kern 1pt 2})\,t^8+(q^{\kern 0.5pt 8}+10\,p^{\kern 1pt 2}\,q^{\kern 0.5pt 6}
+4\,p^{\kern 1pt 4}\,q^{\kern 0.5pt 4}\,-\\
-\,14\,p^{\kern 1pt 6}\,q^{\kern 0.7pt 2}+p^{\kern 1pt 8})\,t^6
-p^{\kern 1pt 2}\,q^{\kern 0.7pt 2}\,(q^{\kern 0.5pt 8}
-14\,p^{\kern 1pt 2}\,q^{\kern 0.5pt 6}+4\,p^{\kern 1pt 4}\,q^{\kern 0.5pt 4}
+10\,p^{\kern 1pt 6}\,q^{\kern 0.7pt 2}+\\
+p^{\kern 1pt 8})\,t^4-p^{\kern 1pt 6}\,q^{\kern 0.5pt 6}\,(q^{\kern 0.7pt 2}
+2\,p^{\kern 1pt 2})\,(3\,p^{\kern 1pt 2}-2\,q^{\kern 0.7pt 2})\,t^2
-q^{\kern 0.7pt 10}\,p^{\kern 1pt 10}=0.
\endgathered
\quad
\mytag{1.15}
$$
\myconjecture{1.1} For any positive coprime integers $p\neq q$ the 
polynomial $Q_{p\kern 0.6pt q}(t)$ in \mythetag{1.12} is irreducible in 
the ring $\Bbb Z[t]$. 
\endproclaim
     This conjecture is known as the second cuboid conjecture. It was 
formulated in \mycite{45}. In particular it means that if it is true, the 
equation \mythetag{1.15} has no integer roots for any positive coprime integers 
$p\neq q$. Like in \mycite{49}, here we shall not try to prove or disprove 
Conjecture~\mytheconjecture{1.1}. Instead, we study real positive roots of the 
equation \mythetag{1.15} in the case where $p$ is much larger than $q$. This
case is reverse to that of \mycite{49}. Below we find asymptotic expansions 
and estimates for the roots of the equation \mythetag{1.15} as $p\to+\infty$.
They can be applied to a numeric search for perfect cuboids. In particular
they can be applied to further improvements of the numeric algorithm 
proposed in \mycite{49}.\par
\head
2. The parameters reversion formula.
\endhead
     In order to obtain asymptotic expansions for roots of the equation
\mythetag{1.15} as $p\to+\infty$ one can go through the same steps as in 
\mycite{49} for the case where $q\to+\infty$. However, this is too expensive 
way. In order to avoid this way below we use the reversion formula derived 
in \mycite{47}. Here is this formula:
$$
Q_{p\kern 0.6pt q}(t)=-\frac{Q_{qp}(p^{\kern 1pt 2}\,q^{\kern 0.7pt 2}/t)\,t^{10}}
{p^{\kern 1pt 10}\,q^{\kern 0.7pt 10}}.
\mytag{2.1}
$$\par
     The formula \mythetag{2.1} relates the polynomial $Q_{p\kern 0.6pt q}(t)$ 
with the reverse polynomial $Q_{qp}(t)$. Both polynomials are even. Therefore, 
like in \mycite{49}, we use the condition 
$$
\hskip -2em
\cases \text{$t>0$ \ if \ $t$ \ is a real root,}\\
\text{$\Rea(t)\geqslant 0$ \ and \ $\Img(t)>0$ \ if \ $t$ \ 
is a complex root}
\endcases
\mytag{2.2}
$$
in order to divide their roots into two groups. Let's denote through
$t_1$, $t_2$, \dots, $t_{10}$ the roots of the polynomial 
$Q_{p\kern 0.6pt q}(t)$ and through $\tilde t_1$, $\tilde t_2$, \dots, 
$\tilde t_{10}$ the roots of the reverse polynomial $Q_{qp}(t)$. Then
the roots $t_1$, $t_2$, \dots, $t_5$ and $\tilde t_1$, $\tilde t_2$, 
\dots, $\tilde t_5$ are those roots that satisfy the condition 
\mythetag{2.2}. Other roots are opposite to them:
 $$
\xalignat 5
&\hskip -1em
t_6=-t_1,&&t_7=-t_2,&&t_8=-t_3,&&t_9=-t_4,&&t_{10}=-t_5,
\qquad\quad\\
\vspace{-1.5ex}
\mytag{2.3}\\
\vspace{-1.5ex}
&\hskip -1em
\tilde t_6=-\tilde t_1,&&\tilde t_7=-\tilde t_2,&&\tilde t_8=
-\tilde t_3,&&\tilde t_9=-\tilde t_4,&&\tilde t_{10}=-\tilde t_5.
\qquad\quad
\endxalignat
$$
Since $q^{\kern 0.7pt 10}\,p^{\kern 1pt 10}\neq 0$ in \mythetag{1.12},
the roots of both polynomials in \mythetag{2.3} are nonzero. Hence the
formula \mythetag{2.1} is applicable to all of them.\par
\head
3. Roots of the reverse polynomial.
\endhead
Note that $p\to+\infty$ for the reverse polynomial $Q_{qp}(t)$ is equivalent 
to $q\to+\infty$ for the original polynomial $Q_{p\kern 0.6pt q}(t)$. Therefore
we can use the results of \mycite{49} and reformulate them as the results
for the roots $\tilde t_1$, \dots, $\tilde t_{10}$ of the
reverse polynomial $Q_{qp}(t)$ as $p\to+\infty$. These results are formulated
using five asymptotic intervals. Three of these five intervals are on the real 
axis: 
$$
\align
\hskip -2em
q^{\kern 0.7pt 2}-\frac{5\,q^{\kern 0.7pt 3}}{p}<\ &\tilde t<q^{\kern 0.7pt 2},
\mytag{3.1}\\
\hskip -2em
q^{\kern 0.7pt 2}<\ &\tilde t<q^{\kern 0.7pt 2}+\frac{5\,q^{\kern 0.7pt 3}}{p}
\mytag{3.2}\\
\hskip -2em
q\,p-\frac{16\,q^{\kern 0.7pt 3}}{p}-\frac{5\,q^{\kern 0.7pt 4}}{p^{\kern 1pt 2}}
<\ &\tilde t<q\,p-\frac{16\,q^{\kern 0.7pt 3}}{p}+\frac{5\,q^{\kern 0.7pt 4}}
{p^{\kern 1pt 2}}.
\mytag{3.3}
\endalign
$$
The other two intervals are on the imaginary axis of the complex plane:
$$
(\sqrt{2}+1)\,p^{\kern 1pt 2}+(\sqrt{2}-2)\,q^{\kern 0.7pt 2}
-\frac{5\,q^{\kern 0.7pt 3}}{p}<\Img\,\tilde t<(\sqrt{2}+1)
\,p^{\kern 1pt 2}+(\sqrt{2}-2)\,q^{\kern 0.7pt 2}
+\frac{5\,q^{\kern 0.7pt 3}}{p}.
\qquad
\mytag{3.4}
$$
\vskip -4ex
$$
\pagebreak
(\sqrt{2}-1)\,p^{\kern 1pt 2}+(\sqrt{2}+2)\,q^{\kern 0.7pt 2}
-\frac{5\,q^{\kern 0.7pt 3}}{p}<\Img\,\tilde t<(\sqrt{2}-1)
\,p^{\kern 1pt 2}+(\sqrt{2}+2)\,q^{\kern 0.7pt 2}
+\frac{5\,q^{\kern 0.7pt 3}}{p}.
\qquad
\mytag{3.5}
$$
Now Theorem~6.1 from \mycite{49} is formulated as follows. 
\mytheorem{3.1} For $p\geqslant 59\,q$ five roots $\tilde t_1$, $\tilde t_2$, 
$\tilde t_3$, $\tilde t_4$, $\tilde t_5$ of the reverse polynomial
$Q_{qp}(t)$ obeying the condition \mythetag{2.2} are simple. They are located 
within five disjoint intervals \mythetag{3.1}, \mythetag{3.2}, \mythetag{3.3}, 
\mythetag{3.4}, \mythetag{3.5}, one per each interval. 
\endproclaim
     Theorem~\mythetheorem{3.1} means that three roots $\tilde t_1$, $\tilde t_2$, 
$\tilde t_3$ of the reverse polynomial $Q_{qp}(t)$ are real. They are arranged 
as follows:
$$
\hskip -2em
\tilde t_1<\tilde t_2<\tilde t_3
\mytag{3.6}
$$
The other two roots $\tilde t_4$ and $\tilde t_5$ are imaginary. They are arranged
so that 
$$
\hskip -2em
\Img\tilde t_4>\Img\tilde t_5.
\mytag{3.7}
$$ 
We are going to preserve the same arrangement \mythetag{3.6} and \mythetag{3.7}
for the roots of the original polynomial $Q_{p\kern 0.6pt q}(t)$, i\.\,e\. we
write
$$
\xalignat 2
&\hskip -1em
t_1<t_2<t_3,
&&\Img t_4>\Img t_5.
\mytag{3.8}
\endxalignat
$$
Then,	applying \mythetag{2.3} and the reversion formula \mythetag{2.1}, from
\mythetag{3.6}, \mythetag{3.7}, and \mythetag{3.8} we derive the following
correspondence for the roots of $Q_{p\kern 0.6pt q}(t)$ and $Q_{qp}(t)$:
$$
\xalignat 5
&t_1=\frac{p^{\kern 1pt 2}\,q^{\kern 0.7pt 2}}{\tilde t_3},
&&t_2=\frac{p^{\kern 1pt 2}\,q^{\kern 0.7pt 2}}{\tilde t_2},
&&t_3=\frac{p^{\kern 1pt 2}\,q^{\kern 0.7pt 2}}{\tilde t_1},
&&t_4=-\frac{p^{\kern 1pt 2}\,q^{\kern 0.7pt 2}}{\tilde t_5},
&&t_5=-\frac{p^{\kern 1pt 2}\,q^{\kern 0.7pt 2}}{\tilde t_4}.
\qquad\quad
\mytag{3.9}
\endxalignat
$$
Below we use the formulas \mythetag{3.9} in order to derive asymptotic 
expansions and estimates for the roots of the original polynomial 
$Q_{p\kern 0.6pt q}(t)$ as $p\to+\infty$.\par
\head
4. Asymptotics of the real roots. 
\endhead
     According to Theorem~\mythetheorem{3.1}, for $p\geqslant 59\,q$
the root $\tilde t_3$ of the reverse polynomial $Q_{qp}(t)$ belongs
to the third asymptotic interval \mythetag{3.3}. Applying the first 
formula \mythetag{3.9}, we derive the following interval for the root
$t_1$ of the polynomial $Q_{p\kern 0.6pt q}(t)$:
$$
\hskip -2em
\frac{p^{\kern 1pt 4}\,q}{p^{\kern 1pt 3}-16\,q^{\kern 0.7pt 2}\,p
+5\,q^{\kern 0.7pt 3}}<t<\frac{p^{\kern 1pt 4}\,q}{p^{\kern 1pt 3}
-16\,q^{\kern 0.7pt 2}\,p-5\,q^{\kern 0.7pt 3}}.
\mytag{4.1}
$$
The left and right sides of the inequalities \mythetag{4.1} have 
the asymptotic expansions
$$
\hskip -2em
\aligned
&\frac{p^{\kern 1pt 4}\,q}{p^{\kern 1pt 3}-16\,q^{\kern 0.7pt 2}\,p
+5\,q^{\kern 0.7pt 3}}=p\,q+\frac{16\,q^{\kern 0.7pt 3}}{p}
-\frac{5\,q^{\kern 0.7pt 4}}{p^{\kern 1pt 2}}+\frac{256\,q^{\kern 0.7pt 5}}
{p^{\kern 1pt 3}}-\frac{160\,q^{\kern 0.7pt 6}}{p^{\kern 1pt 4}}+\dots\,,\\
\vspace{1ex}
&\frac{p^{\kern 1pt 4}\,q}{p^{\kern 1pt 3}-16\,q^{\kern 0.7pt 2}\,p
-5\,q^{\kern 0.7pt 3}}=p\,q+\frac{16\,q^{\kern 0.7pt 3}}{p}
+\frac{5\,q^{\kern 0.7pt 4}}{p^{\kern 1pt 2}}+\frac{256\,q^{\kern 0.7pt 5}}
{p^{\kern 1pt 3}}+\frac{160\,q^{\kern 0.7pt 6}}{p^{\kern 1pt 4}}+\dots
\endaligned
\mytag{4.2}
$$
as $p\to+\infty$. From \mythetag{4.1} and \mythetag{4.2} one can derive 
an asymptotic expansion for $t_1$:
$$
\hskip -2em
t_1=p\,q+\frac{16\,q^{\kern 0.7pt 3}}{p}+R_1(p,q)\text{\ \ as \ }
p\to+\infty. 
\mytag{4.3}
$$	
Like in \mycite{49}, our goal here is to obtain an estimate of the form 
$$
\hskip -2em
|R_1(p,q)|<\frac{C(q)}{p^{\kern 0.7pt 2}}.
\mytag{4.4}
$$
In order to get such an estimate we substitute 
$$
\hskip -2em
t=p\,q+\frac{16\,q^{\kern 0.7pt 3}}{p}+\frac{c}{p^{\kern 1pt 2}}
\mytag{4.5}
$$
into the equation \mythetag{1.15}. Then we perform another substitution 
into the equation obtained as a result of substituting \mythetag{4.5} into 
\mythetag{1.15}:
$$
\hskip -2em
p=\frac{1}{z}.
\mytag{4.6}
$$
Upon two substitutions \mythetag{4.5} and \mythetag{4.6} and upon removing 
denominators the equation \mythetag{1.15} is presented as a polynomial equation 
in the new variables $c$ and $z$. It is a peculiarity of this equation that
it can be written as  
$$
\hskip -2em
f(c,q,z)=-2\,q^{\kern 0.7pt 5}\,c.
\mytag{4.7}
$$
Here $f(c,q,z)$ is a polynomial of three variables given by an explicit formula. 
However, the formula for $f(c,q,z)$ is rather huge. Therefore it is placed to the 
ancillary file \darkred{{\tt strategy\kern -0.5pt\_\kern 1.5pt 
formulas\_\kern 0.5pt 02.txt}} in a machine-readable form.\par
     Let $p\geqslant 59\,q$ and let the parameter $c$ run over the interval
from  $-5\,q^{\kern 0.7pt 4}$ to $5\,q^{\kern 0.7pt 4}$:
$$
\hskip -2em
-5\,q^{\kern 0.7pt 4}<c<5\,q^{\kern 0.7pt 4}.
\mytag{4.8}
$$
From $p\geqslant 59\,q$ and from \mythetag{4.6} we derive the estimate 
$|z|\leqslant 1/59\,q^{\kern 0.7pt -1}$. Using this estimate and using the 
inequalities \mythetag{4.8}, by means of direct calculations one can derive 
the following estimate for the modulus of the function $f(c,q,z)$:
$$
\hskip -2em
|f(c,q,z)|<3\,q^{\kern 0.7pt 9}.
\mytag{4.9}
$$
For fixed $q$ and $z$ the estimate \mythetag{4.9} means that the left hand 
side of the equation \mythetag{4.7} is a continuous function of $c$ whose
values are within the range from $-3\,q^{\kern 0.7pt 9}$ to $3\,q^{\kern 0.7pt 9}$ 
while $c$ runs over the interval \mythetag{4.8}. The right hand side of the 
equation \mythetag{4.7} is also a continuous function of $c$. It decreases from 
$10\,q^{\kern 0.7pt 9}$ to $-10\,q^{\kern 0.7pt 9}$ in the interval \mythetag{4.8}. 
Therefore somewhere in the interval \mythetag{4.8} there is at least one root 
of the polynomial equation \mythetag{4.7}.\par 
     The parameter $c$ is related to the variable $t$ by means of the formula 
\mythetag{4.5}. Therefore the inequalities \mythetag{4.8} for $c$ 
imply the following inequalities for $t$:
$$
\hskip -2em
p\,q+\frac{16\,q^{\kern 0.7pt 3}}{p}-\frac{5\,q^{\kern 0.7pt 4}}
{p^{\kern 1pt 2}}<t<p\,q+\frac{16\,q^{\kern 0.7pt 3}}{p}
+\frac{5\,q^{\kern 0.7pt 4}}{p^{\kern 1pt 2}}.
\mytag{4.10}
$$
This result is worth enough, it is formulated as a theorem.  
\mytheorem{4.1} For each $p\geqslant 59\,q$ there is at least one 
real root of the equation \mythetag{1.15} satisfying the inequalities 
\pagebreak \mythetag{4.10}.
\endproclaim
     Note that the inequalities \mythetag{4.10} are not derived 
from \mythetag{4.1} and \mythetag{4.2}. They follow from \mythetag{4.7}
and \mythetag{4.9} as described above. These inequalities provide 
an estimate of the form \mythetag{4.4} for the remainder term in the
asymptotic expansion \mythetag{4.3}.\par
     Let's proceed to the roots $t_2$ and $t_3$. According to 
Theorem~\mythetheorem{3.1} for $p\geqslant 59\,q$ the roots $\tilde t_2$ 
and $\tilde t_1$ of the reverse polynomial $Q_{qp}(t)$ belongs
to the intervals \mythetag{3.2} and \mythetag{3.1} respectively. Applying 
the second  and the third formulas \mythetag{3.9}, we derive the following 
two intervals for the roots $t_2$ and $t_3$ of the original polynomial 
$Q_{p\kern 0.6pt q}(t)$:     
$$     
\align
\hskip -2em
\frac{p^{\kern 1pt 3}}{p+5\,q}<\ &t<p^{\kern 1pt 2}
\mytag{4.11}\\
\hskip -2em
p^{\kern 1pt 2}<\ &t<\frac{p^{\kern 1pt 3}}{p-5\,q}.
\mytag{4.12}
\endalign
$$
The left and right hand sides of the above inequalities \mythetag{4.11} and 
\mythetag{4.12} have the following asymptotic expansions as $p\to+\infty$:
$$
\hskip -2em
\aligned
&\frac{p^{\kern 1pt 3}}{p+5\,q}=p^{\kern 1pt 2}-5\,q\,p
+25\,q^{\kern 0.7pt 2}-\frac{125\,q^{\kern 0.7pt 3}}{p}+\dots\,,\\
\vspace{1ex}
&\frac{p^{\kern 1pt 3}}{p-5\,q}=p^{\kern 1pt 2}+5\,q\,p
+25\,q^{\kern 0.7pt 2}+\frac{125\,q^{\kern 0.7pt 3}}{p}+\dots\,.
\endaligned
\mytag{4.13}
$$
From \mythetag{4.11}, \mythetag{4.12}, and \mythetag{4.13} one can derive 
asymptotic formulas for $t_2$ and $t_3$:
$$
\xalignat 2
&\hskip -2em
t_2\sim p^{\kern 1pt 2},
&&t_3\sim p^{\kern 1pt 2}
\mytag{4.14}
\endxalignat 
$$
as $p\to+\infty$. Unfortunately the asymptotic formulas \mythetag{4.14} are 
too rough. They need to be refined. The refinement of the first formula
\mythetag{4.14} looks like 
$$
\hskip -2em
t_2=p^{\kern 1pt 2}-2\,q\,p-2\,q^{\kern 0.7pt 2}+R_2(p,q)
\text{\ \ as \ }p\to+\infty.
\mytag{4.15}
$$
Like in \mythetag{4.4}, we need to derive some estimate of the form 
$$
\hskip -2em
|R_2(p,q)|<\frac{C(q)}{p}
\mytag{4.16}
$$
for the remainder term $R_2(p,q)$ in \mythetag{4.15}. For this purpose 
we substitute 
$$
\hskip -2em
t=p^{\kern 1pt 2}-2\,q\,p-2\,q^{\kern 0.7pt 2}+\frac{c}{p}
\mytag{4.17}
$$
into the equations \mythetag{1.15}. Then we replace $p$ with the new variable
$z$ using the substitution \mythetag{4.6}. As a result of two substitutions
\mythetag{4.17} and \mythetag{4.6} upon removing denominators we get a polynomial 
equation in the new variables $c$ and $z$. As it turns out, this polynomial equation 
can be written in the following way:
$$
\hskip -2em
80\,q^{\kern 0.7pt 4}+\varphi(c,q,z)=-16\,q\,c.
\mytag{4.18}
$$
Here $\varphi(c,q,z)$ is a polynomial of three variables given by an explicit 
formula. The formula for $\varphi(c,q,z)$ is rather huge. Therefore it is placed 
to the ancillary file 
\darkred{{\tt strategy\kern -0.5pt\_\kern 1.5pt formulas\_\kern 0.5pt 02.txt}} 
in a machine-readable form.\par
     Let $p\geqslant 59\,q$ and let the parameter $c$ run over the interval
from  $-9\,q^{\kern 0.7pt 3}$ to $0$:
$$
\hskip -2em
-9\,q^{\kern 0.7pt 3}<c<0.
\mytag{4.19}
$$
From $p\geqslant 59\,q$ and from \mythetag{4.6} we derive the estimate 
$|z|\leqslant 1/59\,q^{\kern 0.7pt -1}$. Using this estimate and using the 
inequalities \mythetag{4.19}, by means of direct calculations one can derive 
the following estimate for the modulus of the function $\varphi(c,q,z)$:
$$
\hskip -2em
|\varphi(c,q,z)|<52\,q^{\kern 0.7pt 4}.
\mytag{4.20}
$$
For fixed $q$ and $z$ the estimate \mythetag{4.20} means that the left hand 
side of the equation \mythetag{4.18} is a continuous function of $c$ taking 
its values within the range from $28\,q^{\kern 0.7pt 4}$ to 
$132\,q^{\kern 0.7pt 4}$ while $c$ runs over the interval \mythetag{4.19}. 
The right hand side of the equation \mythetag{4.18} is also a continuous 
function of $c$. It decreases from $144\,q^{\kern 0.7pt 4}$ to $0$ in the 
interval \mythetag{4.19}. Therefore somewhere in the interval \mythetag{4.19} 
there is at least one root of the polynomial equation \mythetag{4.18}.\par 
     The parameter $c$ in \mythetag{4.18} is related to $t$ by means of the 
formula \mythetag{4.17}. Therefore the inequalities \mythetag{4.19} for $c$ 
imply the following inequalities for $t$:
$$
\hskip -2em
p^{\kern 1pt 2}-2\,q\,p-2\,q^{\kern 0.7pt 2}
-\frac{9\,q^{\kern 0.7pt 3}}{p}<t<p^{\kern 1pt 2}-2\,q\,p
-2\,q^{\kern 0.7pt 2}.
\mytag{4.21}
$$
This result is formulated as the following theorem.  
\mytheorem{4.2} For each $p\geqslant 59\,q$ there is at least one 
real root of the equation \mythetag{1.15} satisfying the inequalities 
\mythetag{4.21}.
\endproclaim
     Note that the inequalities \mythetag{4.21} prove the asymptotic
expansion \mythetag{4.15} and provide an estimate of the form
\mythetag{4.16} for the remainder term $R_2(p,q)$ in it.\par 
    The root $t_3$ of the equation \mythetag{1.15} is handled in a
similar way. The refinement of the second asymptotic formula 
\mythetag{4.14} for this root looks like 
$$
\hskip -2em
t_3=p^{\kern 1pt 2}+2\,q\,p-2\,q^{\kern 0.7pt 2}+R_3(p,q)
\text{\ \ as \ }p\to+\infty.
\mytag{4.22}
$$
Its remainder term $R_3(p,q)$ obeys the estimate of the form
$$
\hskip -2em
|R_3(p,q)|<\frac{C(q)}{p}.
\mytag{4.23}
$$
In order to prove the formulas \mythetag{4.22} and \mythetag{4.23}
we substitute 
$$
\hskip -2em
t=p^{\kern 1pt 2}+2\,q\,p-2\,q^{\kern 0.7pt 2}+\frac{c}{p}
\mytag{4.24}
$$
into the equations \mythetag{1.15}. Then we replace $p$ with the new variable
$z$ using the substitution \mythetag{4.6}. As a result of two substitutions
\mythetag{4.24} and \mythetag{4.6} upon removing denominators we get a polynomial 
equation in the new variables $c$ and $z$. As it turns out, this polynomial equation 
can be written in the following way:
$$
\hskip -2em
80\,q^{\kern 0.7pt 4}+\psi(c,q,z)=16\,q\,c.
\mytag{4.25}
$$
Here $\psi(c,q,z)$ is a polynomial of three variables given by an explicit 
formula. The formula for $\psi(c,q,z)$ is rather huge. Therefore it is placed 
to the ancillary file 
\darkred{{\tt strategy\kern -0.5pt\_\kern 1.5pt formulas\_\kern 0.5pt 02.txt}} 
in a machine-readable form.\par
     Let $p\geqslant 59\,q$ and let the parameter $c$ run over the interval
from  $0$ to $9\,q^{\kern 0.7pt 3}$:
$$
\hskip -2em
0<c<9\,q^{\kern 0.7pt 3}.
\mytag{4.26}
$$
From $p\geqslant 59\,q$ and from \mythetag{4.6} we derive the estimate 
$|z|\leqslant 1/59\,q^{\kern 0.7pt -1}$. Using this estimate and using the 
inequalities \mythetag{4.26}, by means of direct calculations one can derive 
the following estimate for the modulus of the function $\psi(c,q,z)$:
$$
\hskip -2em
|\psi(c,q,z)|<52\,q^{\kern 0.7pt 4}.
\mytag{4.27}
$$
For fixed $q$ and $z$ the estimate \mythetag{4.27} means that the left hand 
side of the equation \mythetag{4.25} is a continuous function of $c$ taking 
its values within the range from $28\,q^{\kern 0.7pt 4}$ to 
$132\,q^{\kern 0.7pt 4}$ while $c$ runs over the interval \mythetag{4.26}. 
The right hand side of the equation \mythetag{4.25} is also a continuous 
function of $c$. It increases from $0$ to $144\,q^{\kern 0.7pt 4}$  in the 
interval \mythetag{4.26}. Therefore somewhere in the interval \mythetag{4.26} 
there is at least one root of the polynomial equation \mythetag{4.25}.\par 
     The parameter $c$ is related to the variable $t$ by means of the formula 
\mythetag{4.24}. Therefore the inequalities \mythetag{4.26} for $c$ 
imply the following inequalities for $t$:
$$
\hskip -2em
p^{\kern 1pt 2}+2\,q\,p-2\,q^{\kern 0.7pt 2}
<t<p^{\kern 1pt 2}+2\,q\,p
-2\,q^{\kern 0.7pt 2}+\frac{9\,q^{\kern 0.7pt 3}}{p}.
\mytag{4.28}
$$
This result is formulated as the following theorem.  
\mytheorem{4.3} For each $p\geqslant 59\,q$ there is at least one 
real root of the equation \mythetag{1.15} satisfying the inequalities 
\mythetag{4.28}.
\endproclaim
     Apart from Theorem~\mythetheorem{4.3}, the inequalities \mythetag{4.28}
prove the asymptotic expansion \mythetag{4.22} and the estimate 
\mythetag{4.23} for the remainder term in it.\par
\head
5. Asymptotics of the complex roots. 
\endhead
     According to Theorem~\mythetheorem{3.1}, for $p\geqslant 59\,q$
the root $\tilde t_5$ of the reverse polynomial $Q_{qp}(t)$ belongs
to the fifth asymptotic interval \mythetag{3.5}. Applying the fourth 
formula \mythetag{3.9}, we derive the following interval for the root
$t_4$ of the polynomial $Q_{p\kern 0.6pt q}(t)$:
$$
\hskip -0.2em
\frac{p^{\kern 1pt 3}\,q^{\kern 0.7pt 2}}{(\sqrt{2}-\!1)\,p^{\kern 1pt 3}\!
+\!(\sqrt{2}+\!2)\,p\,q^{\kern 0.7pt 2}\!+\!5\,q^{\kern 0.7pt 3}}\!<\Img\,t
<\!\frac{p^{\kern 1pt 3}\,q^{\kern 0.7pt 2}}{(\sqrt{2}-\!1)\,p^{\kern 1pt 3}\!
+\!(\sqrt{2}+\!2)\,p\,q^{\kern 0.7pt 2}\!-\!5\,q^{\kern 0.7pt 3}}.
\hskip 0.4em
\mytag{5.1}
$$
Like in \mythetag{4.2} and \mythetag{4.13}, using asymptotic expansions
for both sides of the inequalities \mythetag{5.1}, we derive the
following asymptotic expansion for the root $t_4$:
$$
\hskip -2em
t_4=(\sqrt{2}+1)\,\goth i\,q^{\kern 0.7pt 2}+R_4(p,q)
\text{\ \ as \ }p\to+\infty.
\mytag{5.2}
$$
Here $\goth i=\sqrt{-1}$. The expansion \mythetag{5.2} is analogous 
to \mythetag{4.3}, \mythetag{4.15}, and \mythetag{4.22}. Our goal here 
is to derive an estimate for the remainder term $R_4(p,q)$ of the form 
$$
\hskip -2em
|R_4(p,q)|<\frac{C(q)}{p}.
\mytag{5.3}
$$
In order to obtain such an estimate we substitute 
$$
\hskip -2em
t=(\sqrt{2}+1)\,\goth i\,q^{\kern 0.7pt 2}+\frac{\goth i\,c}
{p^{\kern 1pt 2}}
\mytag{5.4}
$$
into the equation \mythetag{1.15}. Then we apply the substitution
\mythetag{4.6} to the equation obtained as a result of substituting 
\mythetag{5.2} into \mythetag{1.15}. Upon applying two substitutions 
\mythetag{5.4} and \mythetag{4.6} and upon removing denominators 
the equation \mythetag{1.15} is written as a polynomial equation 
in the new variables $c$ and $z$. It can be represented as  
$$
\hskip -2em
\eta(c,q,z)=16\,c.
\mytag{5.5}
$$
Here $\eta(c,q,z)$ is a polynomial of three variables given by an explicit 
formula. The formula for $\eta(c,q,z)$ is rather huge. Therefore it is placed 
to the ancillary file 
\darkred{{\tt strategy\kern -0.5pt\_\kern 1.5pt formulas\_\kern 0.5pt 02.txt}} 
in a machine-readable form.\par
     Let $p\geqslant 59\,q$ and let the parameter $c$ run over the interval
from $-5\,q^{\kern 0.7pt 3}$ to $5\,q^{\kern 0.7pt 3}$:
$$
\hskip -2em
-5\,q^{\kern 0.7pt 3}<c<5\,q^{\kern 0.7pt 3}.
\mytag{5.6}
$$
From $p\geqslant 59\,q$ and from \mythetag{4.6} we derive the estimate 
$|z|\leqslant 1/59\,q^{\kern 0.7pt -1}$. Using this estimate and using the 
inequalities \mythetag{5.6}, by means of direct calculations one can derive 
the following estimate for the modulus of the function $\eta(c,q,z)$:
$$
\hskip -2em
|\eta(c,q,z)|<14\,q^{\kern 0.7pt 3}.
\mytag{5.7}
$$
For fixed $q$ and $z$ the estimate \mythetag{5.7} means that the left hand 
side of the equation \mythetag{5.5} is a continuous function of $c$ taking 
its values within the range from $-14\,q^{\kern 0.7pt 3}$ to 
$14\,q^{\kern 0.7pt 3}$ while $c$ runs over the interval \mythetag{5.6}. 
The right hand side of the equation \mythetag{5.5} is also a continuous 
function of $c$. It increases from $-80\,q^{\kern 0.7pt 3}$ to 
$80\,q^{\kern 0.7pt 3}$ in the 
interval \mythetag{5.6}. Therefore somewhere in the interval \mythetag{5.6} 
there is at least one root of the polynomial equation \mythetag{5.5}.\par 
     The parameter $c$ is related to the variable $t$ by means of the formula 
\mythetag{5.4}. Therefore the inequalities \mythetag{5.6} for $c$ 
imply the following inequalities for $t$:
$$
\hskip -2em
(\sqrt{2}+1)\,q^{\kern 0.7pt 2}-\frac{5\,q^{\kern 0.7pt 3}}
{p^{\kern 1pt 2}}<\Img\,t<(\sqrt{2}+1)\,q^{\kern 0.7pt 2}
+\frac{5\,q^{\kern 0.7pt 3}}{p^{\kern 1pt 2}}. 
\mytag{5.8}
$$
This result leads to the following theorem.  
\mytheorem{5.1} For each $p\geqslant 59\,q$ there is at least one 
purely imaginary root of the equation \mythetag{1.15} satisfying 
the inequalities \mythetag{5.8}.
\endproclaim
     The complex root $t_5$ of the polynomial $Q_{p\kern 0.6pt q}(t)$ is 
similar to the root $t_4$. The asymptotic expansion for this root is 
derived from the inequalities \mythetag{3.4} with the use of the fifth
formula \mythetag{3.9}. This expansion looks like 
$$
\hskip -2em
t_5=(\sqrt{2}-1)\,\goth i\,q^{\kern 0.7pt 2}+R_5(p,q)
\text{\ \ as \ }p\to+\infty.
\mytag{5.9}
$$
Here $\goth i=\sqrt{-1}$. The expansion \mythetag{5.9} is analogous 
to \mythetag{5.2}. Our goal here is to derive an estimate for the remainder 
term $R_5(p,q)$ of the form 
$$
\hskip -2em
|R_5(p,q)|<\frac{C(q)}{p}.
\mytag{5.10}
$$
In order to obtain such an estimate we substitute 
$$
\hskip -2em
t=(\sqrt{2}-1)\,\goth i\,q^{\kern 0.7pt 2}+\frac{\goth i\,c}
{p^{\kern 1pt 2}}
\mytag{5.11}
$$
into the equation \mythetag{1.15}. Then we apply the substitution
\mythetag{4.6} to the equation obtained as a result of substituting 
\mythetag{5.11} into \mythetag{1.15}. Upon applying two substitutions 
\mythetag{5.11} and \mythetag{4.6} and upon removing denominators 
the equation \mythetag{1.15} is written as a polynomial equation 
in the variables $c$ and $z$. It can be represented as  
$$
\hskip -2em
\zeta(c,q,z)=16\,c.
\mytag{5.12}
$$
Here $\zeta(c,q,z)$ is a polynomial of three variables given by an explicit 
formula. The formula for $\zeta(c,q,z)$ is rather huge. Therefore it is placed 
to the ancillary file 
\darkred{{\tt strategy\kern -0.5pt\_\kern 1.5pt formulas\_\kern 0.5pt 02.txt}} 
in a machine-readable form.\par
     Let $p\geqslant 59\,q$ and let the parameter $c$ run over the interval
\mythetag{5.6}. From $p\geqslant 59\,q$ and from \mythetag{4.6} we derive 
the estimate 
$|z|\leqslant 1/59\,q^{\kern 0.7pt -1}$. Using this estimate and using the 
inequalities \mythetag{5.6}, by means of direct calculations one can derive 
the following estimate for the modulus of the function $\zeta(c,q,z)$:
$$
\hskip -2em
|\zeta(c,q,z)|<14\,q^{\kern 0.7pt 3}.
\mytag{5.13}
$$
For fixed $q$ and $z$ the estimate \mythetag{5.13} means that the left hand 
side of the equation \mythetag{5.12} is a continuous function of $c$ taking 
its values within the range from $-14\,q^{\kern 0.7pt 3}$ to 
$14\,q^{\kern 0.7pt 3}$ while $c$ runs over the interval \mythetag{5.6}. 
The right hand side of the equation \mythetag{5.12} is also a continuous 
function of $c$. It increases from $-80\,q^{\kern 0.7pt 3}$ to 
$80\,q^{\kern 0.7pt 3}$ in the 
interval \mythetag{5.6}. Therefore somewhere in the interval \mythetag{5.6} 
there is at least one root of the polynomial equation \mythetag{5.12}.\par 
     The parameter $c$ is related to the variable $t$ by means of the formula 
\mythetag{5.11}. Therefore the inequalities \mythetag{5.6} for $c$ 
imply the following inequalities for $t$:
$$
\hskip -2em
(\sqrt{2}-1)\,q^{\kern 0.7pt 2}-\frac{5\,q^{\kern 0.7pt 3}}
{p^{\kern 1pt 2}}<\Img\,t<(\sqrt{2}-1)\,q^{\kern 0.7pt 2}
+\frac{5\,q^{\kern 0.7pt 3}}{p^{\kern 1pt 2}}. 
\mytag{5.14}
$$
This result leads to the following theorem.  
\mytheorem{5.2} For each $p\geqslant 59\,q$ there is at least one 
purely imaginary root of the equation \mythetag{1.15} satisfying 
the inequalities \mythetag{5.14}.
\endproclaim
    Theorems~\mythetheorem{5.1} and \mythetheorem{5.2} solve the problem of 
obtaining estimates of the form \mythetag{5.3} and
\mythetag{5.10} for the remainder terms in the asymptotic expansions
\mythetag{5.2} and \mythetag{5.9} for $p\geqslant 59\,q$. Along with 
Theorems~\mythetheorem{4.1}, \mythetheorem{4.2}, and \mythetheorem{4.3}, 
Theorems~\mythetheorem{5.1} and \mythetheorem{5.2} separate the roots $t_1$, 
$t_2$, $t_3$, $t_4$, $t_5$ of the equation \mythetag{1.15} from each other 
for sufficiently large $p$ and specify their locations.\par
\head
6. Non-intersection of asymptotic intervals.
\endhead
     The roots $\tilde t_1$, $\tilde t_2$, $\tilde t_3$, $\tilde t_4$, 
$\tilde t_5$ of the reverse polynomial belong to the intervals
\mythetag{3.1}, \mythetag{3.2}, \mythetag{3.3}, \mythetag{3.4}, \mythetag{3.5},
one per each interval. However, the asymptotic intervals \mythetag{4.10}, 
\mythetag{4.21}, \mythetag{4.28}, \mythetag{5.8}, \mythetag{5.14} do not 
exactly correspond to them by virtue of the formula \mythetag{3.9}. For this 
reason we need to prove some non-intersection results concerning the intervals 
\mythetag{4.10}, \mythetag{4.21}, \mythetag{4.28}, \mythetag{5.8}, and
\mythetag{5.14}.
\mylemma{6.1} For $p\geqslant 59\,q$ the asymptotic intervals 
\mythetag{4.10}, \mythetag{4.21}, \mythetag{4.28}, \mythetag{5.8}, and
\mythetag{5.14} do not comprise the origin. 
\endproclaim
\demo{Proof} Indeed, from $p\geqslant 59\,q$ for the left endpoint of the 
interval \mythetag{4.10} we derive
$$
\hskip -2em
p\,q+\frac{16\,q^{\kern 0.7pt 3}}{p}-\frac{5\,q^{\kern 0.7pt 4}}
{p^{\kern 1pt 2}}>p\,q-\frac{5\,q^{\kern 0.7pt 4}}{p^{\kern 1pt 2}}
\geqslant 59\,q^{\kern 0.7pt 2}-\frac{5\,q^{\kern 0.7pt 2}}
{59^{\kern 0.7pt 2}}>58\,q^{\kern 0.7pt 2}>0.
\mytag{6.1}
$$
For the left endpoint of the interval \mythetag{4.21} we derive
$$
\gathered
p^{\kern 1pt 2}-2\,q\,p-2\,q^{\kern 0.7pt 2}
-\frac{9\,q^{\kern 0.7pt 3}}{p}=(p-q)^2-3\,q^{\kern 0.7pt 2}
-\frac{9\,q^{\kern 0.7pt 3}}{p}\geqslant\\
\geqslant (58\,q)^2-3\,q^{\kern 0.7pt 2}-\frac{9\,q^{\kern 0.7pt 2}}{59}
>3360\,q^{\kern 0.7pt 2}>0.
\endgathered
\mytag{6.2}
$$
The case of the interval \mythetag{4.28} is similar. In this case we have 
$$
\hskip -2em
p^{\kern 1pt 2}+2\,q\,p-2\,q^{\kern 0.7pt 2}=
(p+q)^2-3\,q^{\kern 0.7pt 2}\geqslant(60\,q)^2-3\,q^{\kern 0.7pt 2}=
3597\,q^{\kern 0.7pt 2}>0.
\mytag{6.3}
$$
For the bottom endpoints of the intervals 
\mythetag{5.8} and \mythetag{5.14} from $p\geqslant 59\,q$ we derive 
$$
\xalignat 2
&(\sqrt{2}+1)\,q^{\kern 0.7pt 2}
-\frac{5\,q^{\kern 0.7pt 3}}
{p^{\kern 1pt 2}}>2.4\,q^{\kern 0.7pt 2}>0,
&&(\sqrt{2}-1)\,q^{\kern 0.7pt 2}
-\frac{5\,q^{\kern 0.7pt 3}}
{p^{\kern 1pt 2}}>0.4\,q^{\kern 0.7pt 2}>0.
\quad
\mytag{6.4}
\endxalignat
$$
The above inequalities \mythetag{6.1}, \mythetag{6.2}, \mythetag{6.3}, and
\mythetag{6.4} prove Lemma~\mythelemma{6.1}.
\qed\enddemo
     Lemma~\mythelemma{6.1} means that for $p\geqslant 59\,q$ the real 
intervals \mythetag{4.10}, \mythetag{4.21}, and \mythetag{4.28} do not 
intersect with the imaginary intervals \mythetag{5.8} and \mythetag{5.14}. 
Moreover, the inequalities \mythetag{6.1}, \mythetag{6.2}, \mythetag{6.3}
and \mythetag{6.4} show that all of these intervals are located within 
positive half-lines of the real and imaginary axes. Therefore any roots of 
the equation  \mythetag{1.15} enclosed within these intervals satisfy
the condition \mythetag{2.2}. 
\mylemma{6.2} For $p\geqslant 59\,q$ the asymptotic intervals 
\mythetag{4.10}, \mythetag{4.21}, \mythetag{4.28}, \mythetag{5.8},
and \mythetag{5.14}, do not intersect with each other. 
\endproclaim
\demo{Proof} Let's compare the left endpoint of the interval \mythetag{4.21}
with the right endpoint of the interval \mythetag{4.10}. For their difference
we have the inequalities 
$$
\gather
\biggl(p^{\kern 1pt 2}-2\,q\,p-2\,q^{\kern 0.7pt 2}
-\frac{9\,q^{\kern 0.7pt 3}}{p}\biggr)
-\biggl(p\,q+\frac{16\,q^{\kern 0.7pt 3}}{p}
+\frac{5\,q^{\kern 0.7pt 4}}{p^{\kern 1pt 2}}\biggr)=\\
=\Bigl(p-\frac{3\,q}{2}\Bigr)^2-\frac{17\,q^{\kern 0.7pt 2}}{4}
-\frac{25\,q^{\kern 0.7pt 3}}{p}
-\frac{5\,q^{\kern 0.7pt 4}}{p^{\kern 1pt 2}}\geqslant
\mytag{6.5}\\
\vspace{0.1ex}
\geqslant\Bigl(\frac{115\,q}{2}\Bigr)^2
-\frac{17\,q^{\kern 0.7pt 2}}{4}-\frac{25\,q^{\kern 0.7pt 2}}{59}
-\frac{5\,q^{\kern 0.7pt 2}}{59^{\kern 0.7pt 2}}
>3301\,q^{\kern 0.7pt 2}>0.
\endgather
$$
Similarly, let's compare the left endpoint of the interval \mythetag{4.28}
with the right endpoint of the interval \mythetag{4.21}. For their difference
we have the inequalities 
$$
\hskip -2em
\bigl(p^{\kern 1pt 2}+2\,q\,p-2\,q^{\kern 0.7pt 2}\bigr)
-\bigl(p^{\kern 1pt 2}-2\,q\,p-2\,q^{\kern 0.7pt 2}\bigr)
=4\,q\,p\geqslant4\cdot 59\,q^{\kern 0.7pt 2}>0.
\mytag{6.6}
$$
In the case of imaginary intervals we compare the bottom endpoint of 
\pagebreak the interval \mythetag{5.8} with the top endpoint of the 
interval \mythetag{5.14}. For their difference we have 
$$
\hskip -2em
\gathered
\biggl((\sqrt{2}+1)\,q^{\kern 0.7pt 2}-\frac{5\,q^{\kern 0.7pt 3}}
{p^{\kern 1pt 2}}\biggr)-\biggl((\sqrt{2}-1)\,q^{\kern 0.7pt 2}
+\frac{5\,q^{\kern 0.7pt 3}}{p^{\kern 1pt 2}}\biggr)=\\
=\,2\,q^{\kern 0.7pt 2}-\frac{10\,q^{\kern 0.7pt 3}}
{p^{\kern 1pt 2}}\geqslant 2\,q^{\kern 0.7pt 2}-
\frac{10\,q^{\kern 0.7pt 2}}{59^{\kern 0.7pt 2}}
>q^{\kern 0.7pt 2}>0.
\endgathered
\mytag{6.7}
$$
The above inequalities \mythetag{6.5}, \mythetag{6.6}, and
\mythetag{6.7} prove Lemma~\mythelemma{6.2}.
\qed\enddemo
     Lemmas~\mythelemma{6.1} and \mythelemma{6.2} are summed up in the 
following theorem.
\mytheorem{6.1} For $p\geqslant 59\,q$ five roots $t_1$, $t_2$, $t_3$, 
$t_4$, $t_5$ of the equation \mythetag{1.15} obeying the condition 
\mythetag{2.2} are simple. They are located within five disjoint intervals 
\mythetag{4.10}, \mythetag{4.21}, \mythetag{4.28}, \mythetag{5.8},
\mythetag{5.14}, one per each interval. 
\endproclaim
     Due to \mythetag{2.3} Theorem~\mythetheorem{6.1} locates all of the ten
roots of the equation \mythetag{1.15}. 
\head
7. Integer points of asymptotic intervals.
\endhead
     The intervals \mythetag{4.10}, \mythetag{4.21}, \mythetag{4.28}, 
\mythetag{5.8}, and \mythetag{5.14} are asymptotically small. Their 
length decrease as $p\to+\infty$. Like in \mycite{49}, we use this fact 
in order to determine the number of integer points within them. 
\mytheorem{7.1} If $p\geqslant 59\,q$ and $p>9\,q^{\kern 0.7pt 3}$, then 
the asymptotic intervals \mythetag{4.21} and \mythetag{4.28} have no 
integer points. 
\endproclaim
\mytheorem{7.2} If $p\geqslant 59\,q$ and $p^{\kern 1pt 2}
>10\,q^{\kern 0.7pt 4}$, then the asymptotic interval \mythetag{4.10} 
has at most one integer point. 
\endproclaim
Theorems~\mythetheorem{7.1} and \mythetheorem{7.2} are immediate from the 
inequalities $9\,q^{\kern 0.7pt 3}/p<1$ and 
$10\,q^{\kern 0.7pt 4}/p^{\kern 1pt 2}<1$ that follow from 
$p>9\,q^{\kern 0.7pt 3}$ and $p^{\kern 1pt 2}>10\,q^{\kern 0.7pt 4}$
respectively. We preserve the inequality $p\geqslant 59\,q$ in
Theorems~\mythetheorem{7.1} and \mythetheorem{7.2} in order to emphasize
their relation to the roots of the equation \mythetag{1.15} through 
Theorem~\mythetheorem{6.1}.\par
\mytheorem{7.3} If $p\geqslant 59\,q$ and $p\geqslant 16\,q^{\kern 0.7pt 3}
+5\,q/16$, then the asymptotic interval \mythetag{4.10} has no integer points. 
\endproclaim
     Theorem~\mythetheorem{7.3} is more complicated than 
Theorems~\mythetheorem{7.1} and \mythetheorem{7.2}. But its proof repeats
the arguments used in proving Theorem~7.3 in \mycite{49}. For this reason we
do not provide its proof here. 
\head
8. Application to the cuboid problem.
\endhead
     The equation \mythetag{1.15} is a reduced version of the cuboid
characteristic equation \mythetag{1.8}. It is related to the perfect cuboid
problem through Theorem~\mythetheorem{1.1}. 
Substituting either \mythetag{1.10} or \mythetag{1.13} into the inequalities 
$t>a$, $t>b$, and $t>u$ from Theorem~\mythetheorem{1.1}, we get the following
result expressed by the inequalities 
$$
\xalignat 3
&\hskip -2em
t>p^{\kern 1pt 2},
&&t>p\,q,
&&t>q^{\kern 0.7pt 2}.
\mytag{8.1}
\endxalignat
$$
Similarly, substituting either \mythetag{1.10} or \mythetag{1.13} 
into the inequality $(a+t)\,(b+t)>2\,t^2$ from Theorem~\mythetheorem{1.1}, 
we get the result expressed by the inequality
$$
\pagebreak
\hskip -2em
(p^{\kern 1pt 2}+t)\,(p\,q+t)>2\,t^2.
\mytag{8.2}
$$
Theorem~\mythetheorem{1.1} specified for the case of second cuboid conjecture
(see Conjecture~\mytheconjecture{1.1}) is formulated in the following way.
\mytheorem{8.1} A triple of integer numbers $p$, $q$, and $t$ satisfying the
equation \mythetag{1.15} and such that $p\neq q$ are coprime provides a perfect
cuboid if and only if the inequalities \mythetag{8.1} and \mythetag{8.2} are
fulfilled. 
\endproclaim
     The inequalities \mythetag{8.1} set lower bounds for $t$. The inequality 
\mythetag{8.2} is different. It sets the upper bound for $t$ in the form
of the irrational inequality 
$$
t<\frac{p^{\kern 1pt 2}+p\,q}{2}+\frac{p\,\sqrt{p^{\kern 1pt 2}+6\,p\,q
+q^{\kern 0.7pt 2}}}{2}.
$$\par
     Assume that the inequality $p\geqslant 59\,q$ is fulfilled and assume that
$t$ belongs to the second asymptotic interval \mythetag{4.21}. The inequality
 $t<p^{\kern 1pt 2}-2\,q\,p-2\,q^{\kern 0.7pt 2}$ from \mythetag{4.21} and 
 the inequality $t>p^{\kern 1pt 2}$ from \mythetag{8.1} imply the inequality 
$p^{\kern 1pt 2}<p^{\kern 1pt 2}-2\,q\,p-2\,q^{\kern 0.7pt 2}$. This inequality is 
contradictory since $2\,q\,p+2\,q^{\kern 0.7pt 2}>0$. The contradiction obtained
proves the following theorem.
\mytheorem{8.2} If $p\geqslant 59\,q$, then the asymptotic interval \mythetag{4.21} 
has no points satisfying the inequalities \mythetag{8.1}.
\endproclaim
     Now assume that the inequality $p\geqslant 59\,q$ is fulfilled and assume 
that $t$ belongs to the first asymptotic interval \mythetag{4.10}. In this case 
we have the following two inequalities taken from \mythetag{8.1} and 
\mythetag{4.10} respectively:
$$
\xalignat 2
&\hskip -2em
p^{\kern 1pt 2}<t,
&&t<p\,q+\frac{16\,q^{\kern 0.7pt 3}}{p}
+\frac{5\,q^{\kern 0.7pt 4}}{p^{\kern 1pt 2}}.
\mytag{8.3}
\endxalignat
$$
The inequalities \mythetag{8.3} imply an inequality for $p$ and $q$ without $t$:
$$
\hskip -2em
p^{\kern 1pt 2}<p\,q+\frac{16\,q^{\kern 0.7pt 3}}{p}
+\frac{5\,q^{\kern 0.7pt 4}}{p^{\kern 1pt 2}}.
\mytag{8.4}
$$
The inequality \mythetag{8.4} can be transformed in the following way:
$$
\hskip -2em
\Bigl(p-\frac{q}{2}\Bigr)^2-\frac{q^{\kern 0.7pt 2}}{4}
-\frac{16\,q^{\kern 0.7pt 3}}{p}
-\frac{5\,q^{\kern 0.7pt 4}}{p^{\kern 1pt 2}}<0.
\mytag{8.5}
$$
On the other hand, applying $p\geqslant 59\,q$ to the left hand side of
\mythetag{8.5}, we derive 
$$
\hskip -2em
\Bigl(p-\frac{q}{2}\Bigr)^2-\frac{q^{\kern 0.7pt 2}}{4}
-\frac{16\,q^{\kern 0.7pt 3}}{p}
-\frac{5\,q^{\kern 0.7pt 4}}{p^{\kern 1pt 2}}
\geqslant 3421\,q^{\kern 0.7pt 2}>0.
\mytag{8.6}
$$
The inequalities \mythetag{8.5} and \mythetag{8.6} contradict each other.
The contradiction obtained proves the following theorem. 
\mytheorem{8.3} If $p\geqslant 59\,q$, then the asymptotic interval 
\mythetag{4.10} has no points satisfying the inequalities \mythetag{8.1}.
\endproclaim
     Theorem~\mythetheorem{8.3} complements Theorem~\mythetheorem{7.3}. 
Similarly, Theorem~\mythetheorem{8.2} complements Theorem~\mythetheorem{7.2}
in the case of the asymptotic interval \mythetag{4.21}. \pagebreak
These theorems along with Theorem~\mythetheorem{8.1} are summarized in 
the following theorem. \par
\mytheorem{8.4} If $p\geqslant 59\,q$, then the Diophantine equation 
\mythetag{1.15} has no solutions providing perfect cuboids outside 
the third asymptotic interval \mythetag{4.28}. 
\endproclaim
     The third asymptotic interval \mythetag{4.28} is exceptional in 
Theorem~\mythetheorem{8.4}. The inequality $p\geqslant 59\,q$ does not 
contradict the inequalities \mythetag{8.1} and \mythetag{8.2} within
this interval. However, this interval is cut off by applying 
Theorem~\mythetheorem{7.1} to it. This yields the following result. 
\mytheorem{8.5} If $p\geqslant 59\,q$ and $p>9\,q^{\kern 0.7pt 3}$, 
then the Diophantine equation \mythetag{1.15} has no solutions 
providing perfect cuboids at all. 
\endproclaim
\head
9. Conclusions. 
\endhead
     Theorems~\mythetheorem{8.4} and ~\mythetheorem{8.5} constitute
the main result of the present paper. They should be complemented
with Theorem~8.5 from \mycite{49}. It is formulated as follows.  
\mytheorem{9.1} If $q\geqslant 59\,p$, then the Diophantine equation 
\mythetag{1.15} has no solutions providing perfect cuboids. 
\endproclaim
Theorems~\mythetheorem{8.4} and ~\mythetheorem{8.5} along with 
Theorem~\mythetheorem{9.1} proved in \mycite{49} outline three regions 
in the positive quadrant of the $p\,q$\,-\,coordinate plane. These 
regions are:
\roster
\item"1)" {\bf linear region} given by the linear inequalities
$$
\xalignat 2
&\frac{q}{59}<p,
&&p<59\,q;
\mytag{9.1}
\endxalignat
$$
\item"2)" {\bf nonlinear region} given by the nonlinear inequalities
$$
\xalignat 2
&59\,q\leqslant p,
&&p\leqslant 9\,q^{\kern 0.7pt 3};
\mytag{9.2}
\endxalignat
$$
\item"3)" {\bf no cuboid region} which is the rest of the positive 
$p\,q$\,-\,quadrant.
\endroster
The inequalities \mythetag{9.1} and \mythetag{9.2}, as well as other more
special inequalities of this paper, could be used in order to further 
optimize algorithms for a numeric search of perfect cuboids in the case 
of the second cuboid conjecture. 

\enddocument
\end